\documentclass[preprint,3p,times]{elsarticle}

\usepackage{amssymb}
\usepackage{amsmath}
\usepackage{amsthm}
\usepackage{color}
\usepackage{enumitem,verbatim}
\usepackage{tikz}
\usepackage{appendix}
\usepackage[colorlinks=true]{hyperref}

\theoremstyle{plain}
  \newtheorem{thm}{Theorem}[section]
  \newtheorem{lem}[thm]{Lemma}
  \newtheorem{prop}[thm]{Proposition}
  \newtheorem{cor}[thm]{Corollary}
\theoremstyle{definition}
  
  \newtheorem{exmp}[thm]{Example}

\DeclareMathAlphabet{\mathcal}{OMS}{cmsy}{m}{n}

\DeclareMathOperator{\Mod}{Mod}
\DeclareMathOperator{\LO}{LO}

\DeclareMathOperator{\id}{id}

\makeatletter
\def\ps@pprintTitle{%
 \let\@oddhead\@empty
 \let\@evenhead\@empty
 \def\@oddfoot{\centerline{\thepage}}%
 \let\@evenfoot\@oddfoot}
\makeatother

\renewcommand{\phi}{\varphi}

\newcommand{\CC}{\mathcal{C}}
\newcommand{\CD}{\mathcal{D}}
\newcommand{\CI}{\mathcal{I}}
\newcommand{\CK}{\mathcal{K}}

\newcommand{\CS}{\mathcal{S}}
\newcommand{\CT}{\mathcal{T}}
\newcommand{\CU}{\mathcal{U}}
\newcommand{\CV}{\mathcal{V}}
\newcommand{\CW}{\mathcal{W}}

\newcommand{\bbQ}{\mathbb{Q}}
\newcommand{\bbR}{\mathbb{R}}

\newcommand{\FS}{\mathfrak{S}}

\newcommand{\with}{\mathrel{\&}}

\renewcommand{\leq}{\leqslant}
\renewcommand{\geq}{\geqslant}

\numberwithin{equation}{section}

\allowdisplaybreaks

\begin{document}

\begin{frontmatter}



\title{The complexity of classifying continuous t-norms up to isomorphism}


\author{Jialiang He}
\ead{jialianghe@scu.edu.cn}

\author{Lili Shen}
\ead{shenlili@scu.edu.cn}

\author{Yi Zhou\corref{cor}}
\ead{zhouyi.math@outlook.com}

\cortext[cor]{Corresponding author.}
\address{School of Mathematics, Sichuan University, Chengdu 610064, China}

\begin{abstract}
It is shown that the isomorphism relation between continuous t-norms is Borel bireducible with the relation of order isomorphism between linear orders on the set of natural numbers, and therefore, it is a Borel complete equivalence relation.
\end{abstract}


\end{frontmatter}

\section{Introduction}

Triangular norms, usually referred to as \emph{t-norms} \cite{Klement2000,Alsina2006}, are binary operations $*$ on the unit interval $[0,1]$ such that $([0,1],*,1)$ is a commutative ordered monoid. The origins of t-norms can be traced back to the study of \emph{probabilistic metric spaces} \cite{Menger1942,Schweizer1983}, where they were used to generalize the triangle inequality of classical metric spaces. Later, t-norms received wide attention in many-valued logic (see, e.g., \cite{Hajek1998,Metcalfe2007,Cintula2009,Esteva2009}), where they provide a mathematical framework for generalizing the concept of conjunction in logic. Beyond logic, t-norms also serve as an indispensable tool in various fields, including game theory \cite{Butnariu1993}, decision making \cite{Fodor1994,Grabisch1995}, non-additive measures \cite{Pap1995} and statistics \cite{Nelsen2006}.

A t-norm $*$ is called \emph{continuous} if it is continuous as a function with respect to the usual topology on $[0,1]\times[0,1]$. As the most important class of t-norms, continuous t-norms are a primary focus of research in this field. In particular, the \emph{basic fuzzy logic} (or \emph{BL} for short), developed by H{\'a}jek \cite{Hajek1998}, is known as a logic of \emph{continuous t-norms} \cite{Hajek1998,Cignoli2000}. It is well known that every continuous t-norm is an \emph{ordinal sum} of the three basic t-norms \cite{Faucett1955,Mostert1957,Klement2000,Klement2004b,Alsina2006}, i.e, the minimum, the product, and the {\L}ukasiewicz t-norms; we record an equivalent form of this result as Lemma \ref{t-norm-rep}.


Classifying a class of mathematical objects up to isomorphism is a fundamental problem in various areas of mathematics. For example, finite-dimensional vector spaces over a given field 
are classified up to isomorphism by their dimension, and finite fields are classified up to isomorphism by their order. As the structure of a continuous t-norm is clear by decomposing it into an ordinal sum of the three basic t-norms, it is natural to ask whether it is possible to classify continuous t-norms up to isomorphism. A little surprisingly, as revealed in Theorem \ref{t-norm-iso}, as well as in Examples \ref{limit} and \ref{Cantor}, it can be rather challenging to determine whether given continuous t-norms are isomorphic. Hence, we settle for investigating the \emph{complexity} of classifying continuous t-norms up to isomorphism through \emph{descriptive set theory} \cite{Kechris1995,Gao2008}.

Explicitly, the set $\CT$ of continuous t-norms becomes a Polish space by taking the sup metric (Proposition \ref{T-Polish}). We prove the following results:
\begin{itemize}
\item The isomorphism relation $\cong_t$ on $\CT$ is Borel reducible to the isomorphism relation on $\Mod(L_1)$, where $L_1$ is a finite relational language consisting of a binary relation symbol and three unary relation symbols (Proposition \ref{leqB-ModL1}).
\item The equivalence relation $\cong_{L_0}$ of order isomorphism on the set $\LO$ of linear orders on the set of natural numbers is Borel reducible to the isomorphism relation $\cong_t$ on $\CT$ (Proposition \ref{L-leqC}).
\end{itemize}
Our main result, Theorem \ref{main}, then arises from the combination the these two propositions:
\begin{itemize}
\item The isomorphism relation $\cong_t$ on $\CT$ is Borel bireducible with the equivalence relation $\cong_{L_0}$ on $\LO$.
\end{itemize}
As immediate consequences, we deduce that the isomorphism relation $\cong_t$ is a complete analytic subset of $\CT\times\CT$ (Corollary \ref{cong-t-complete-analytic}), and thus it is not smooth (Corollary \ref{congt-not-smooth}). Moreover, it is Borel complete in the sense that it is Borel bireducible with the universal $S_{\infty}$-orbit equivalence relation (Corollary \ref{congt-Borel-complete}).

The contents of this paper are organized as follows. We discuss the isomorphisms of continuous t-norms as preparations in Section \ref{Continuous-t-norms}, and recall the basic facts about the reducibility of equivalence relations in Section \ref{Reducibility}. Subsequently, Section \ref{Complexity} is devoted to investigating the complexity of the isomorphism relation $\cong_t$ on the set $\CT$ of continuous t-norms.


\section{Isomorphisms of continuous t-norms} \label{Continuous-t-norms}

Given an interval $[a,b]\subseteq\bbR$, a continuous function $*\colon[a,b]\times[a,b]\to[a,b]$ is called a \emph{continuous t-norm} on $[a,b]$ \cite{Klement2000,Klement2004b,Alsina2006} if
\begin{itemize}
\item $([a,b],*,b)$ is a commutative monoid, and
\item $p*q\leq p'*q'$ if $p\leq p'$ and $q\leq q'$ in $[a,b]$.
\end{itemize}

For each continuous t-norm $([a,b],*)$ and $q\in[a,b]$, we say that
\begin{itemize}
\item $q$ is \emph{idempotent}, if $q*q=q$;
\item $q$ is \emph{nilpotent}, if $q_*^{(n)}:=\underbrace{q*q*\dots*q}_{n\ \text{times}}=a$ for some $n\in\omega\setminus\{0\}$, where $\omega=\{0,1,2,\dots\}$ refers to the set of natural numbers.
\end{itemize}

\begin{lem} (See \cite[Proposition 2.3]{Klement2000}.) \label{idempotent}
Let $([a,b],*)$ be a continuous t-norm. If $q\in[a,b]$ is idempotent, then $p*q=\min\{p,q\}$ for all $p\in(a,b)$.
\end{lem}

\begin{exmp} \label{*M*L*P-def}
The following continuous t-norms on the unit interval $[0,1]$ are the most prominent ones:
\begin{itemize}
\item The \emph{minimum t-norm} $([0,1],*_M)$ with $p*_M q=\min\{p,q\}$ for all $p,q\in[0,1]$, in which every $q\in(0,1)$ is idempotent and non-nilpotent.
\item The \emph{product t-norm} $([0,1],*_P)$ with $p*_P q=pq$ being the usual product of $p,q\in[0,1]$, in which every $q\in(0,1)$ is non-idempotent and non-nilpotent.
\item The \emph{{\L}ukasiewicz t-norm} $([0,1],*_{\text{\L}})$ with $p*_{\text{\L}} q=\max\{0,p+q-1\}$ for all $p,q\in[0,1]$, in which every $q\in(0,1)$ is non-idempotent and nilpotent.
\end{itemize}
\end{exmp}


We say that continuous t-norms $([a_1,b_1],*_1)$ and $([a_2,b_2],*_2)$ are \emph{isomorphic}, denoted by 
\[([a_1,b_1],*_1)\cong_t([a_2,b_2],*_2),\]
if there exists an order isomorphism 
\[\phi\colon[a_1,b_1]\to[a_2,b_2]\] 
such that
\[\phi\colon([a_1,b_1],*_1,b_1)\to([a_2,b_2],*_2,b_2)\]
is an isomorphism of monoids. It is straightforward to verify the following two lemmas:

\begin{lem}(See \cite[Theorem 2.6]{Klement2004b}.) \label{*-iso-PL-condition:L}
A continuous t-norm $([a,b],*)$ is isomorphic to the {\L}ukasiewicz t-norm $([0,1],*_{\text{\L}})$ if, and only if, every $q\in(a,b)$ is nilpotent.
\end{lem}

\begin{lem} \label{*-iso-PL-condition:min}
Let $([a,b],*)$ be a continuous t-norm. If $([a,b],*)$ is isomorphic to the product t-norm $([0,1],*_P)$ or the {\L}ukasiewicz t-norm $([0,1],*_{\text{\L}})$, then $p*q<\min\{p,q\}$ for all $p,q\in(a,b)$.
\end{lem}




It is well known \cite{Faucett1955,Mostert1957,Klement2000,Klement2004b,Alsina2006} that every continuous t-norm $*$ on the unit interval $[0,1]$ can be written as an \emph{ordinal sum} of the minimum, the product, and the {\L}ukasiewicz t-norm. Explicitly:

\begin{lem} \label{t-norm-rep} {\rm\cite{Klement2000,Klement2004b,Alsina2006}}
For each continuous t-norm $([0,1],*)$, there exists a countable index set $A$ and a family of pairwise disjoint nonempty open intervals $\{(a_{\alpha},b_{\alpha}) : \alpha\in A\}$ such that:
\begin{itemize}
    \item the set of non-idempotent elements is exactly $\bigcup\limits_{\alpha\in A}(a_{\alpha},b_{\alpha})$;
    \item for each $\alpha\in A$, the continuous t-norm $([a_{\alpha},b_{\alpha}],*)$ obtained by restricting $*$ to $[a_{\alpha},b_{\alpha}]$ is either isomorphic to the product t-norm $([0,1],*_P)$ or isomorphic to the {\L}ukasiewicz t-norm $([0,1],*_{\text{\L}})$;
    \item for any $x, y \in [0,1]$, if there is no $\alpha \in A$ such that $\{x, y\} \subset [a_\alpha, b_\alpha]$, then $x * y = \min(x, y)$.
\end{itemize}
In particular, if $([0,1],*)$ is the minimum t-norm $([0,1],*_M)$, then $A=\varnothing$.
\end{lem}

Therefore, it makes sense to denote a continuous t-norm $([0,1],*)$ by 
\begin{equation} \label{t-norm-def-by-rep}
([a_{\alpha}, b_{\alpha}],*)_{\alpha\in A},
\end{equation}
where the intervals $[a_{\alpha},b_{\alpha}]$ $(\alpha\in A)$ are obtained by Lemma \ref{t-norm-rep}. Specifically, the minimum t-norm $([0,1],*_M)$ is denoted by $([a_{\alpha}, b_{\alpha}],*_M)_{\alpha\in \varnothing}$. In order to characterize isomorphisms of continuous t-norms in terms of their interval structures, we define:
\begin{itemize}
\item $\CI_P^*=\{(a_{\alpha},b_{\alpha})\mid \alpha\in A\ \text{and}\ ([a_{\alpha},b_{\alpha}],*)\cong_t([0,1],*_P)\}$.
\item $\CI_{\text{\L}}^*=\{(a_{\alpha},b_{\alpha})\mid \alpha\in A\ \text{and}\ ([a_{\alpha},b_{\alpha}],*)\cong_t([0,1],*_{\text{\L}})\}$.
\item $\CI_M^* = \{(a,b)\mid (a,b)\ \text{is a maximal nonempty open interval of idempotent elements of}\ ([0,1],*)\}$; in other words, $(a,b)\in\CI_M$ if 
    \begin{itemize}
    \item every $q\in(a,b)$ is an idempotent element of $*$, and
    \item for each $\epsilon>0$ with $a-\epsilon\geq 0$, there exist non-idempotent elements of $*$ in $(a-\epsilon,a)$, and for each $\epsilon>0$ with $b+\epsilon\leq 1$, there exist non-idempotent elements of $*$ in $(b,b+\epsilon)$.
    \end{itemize}
\end{itemize}

Let $\FS$ denote the set of all collections of nonempty, pairwise disjoint open subintervals of $[0,1]$ whose union is dense in $[0,1]$. In other words, $\CS\in\FS$ if $\CS$ is a collection of nonempty, pairwise disjoint open subintervals of $[0,1]$ and 
\[\overline{\bigcup\CS}=[0,1].\]
Then:
\begin{itemize}
\item each $\CS\in\FS$ is equipped with a (strict) linear order $\prec$ given by
\begin{equation} \label{CS-order}
(a,b)\prec(c,d)\iff b\leq c
\end{equation}
for all $(a,b),(c,d)\in\CS$;
\item there is a map $\varUpsilon$ from the set $\CT$ of continuous t-norms on $[0,1]$ to $\FS$ given by
\begin{equation} \label{Up-def}
\varUpsilon\colon\CT\to\FS,\quad([0,1],*)=([a_{\alpha},b_{\alpha}],*)_{\alpha\in A}\mapsto\CI_P^*\cup\CI_{\text{\L}}^*\cup\CI_M^*.
\end{equation}
Note that $\varUpsilon([0,1],*)=\CI_P^*\cup\CI_{\text{\L}}^*\cup\CI_M^*$ is never empty by Lemma \ref{t-norm-rep}. In particular, for the minimum t-norm $([0,1],*_M)$, we have $\CI_P^{*_M}=\CI_{\text{\L}}^{*_M}=\varnothing$ and $\CI_M^{*_M}=\{(0,1)\}$, and consequently $\varUpsilon([0,1],*_M)=\{(0,1)\}$.
\end{itemize}

We are now able to characterize isomorphisms of continuous t-norms as follows:

\begin{thm} \label{t-norm-iso}
Let $([0,1],*)=([a_{\alpha},b_{\alpha}],*)_{\alpha\in A}$ and $([0,1],\bullet)=([c_{\beta},d_{\beta}],\bullet)_{\beta\in B}$ be continuous t-norms. Then the following statements are equivalent:
\begin{enumerate}[label={\rm(\arabic*)}]
\item \label{t-norm-iso:cong}
$([0,1],*)\cong_t([0,1],\bullet)$.
\item \label{t-norm-iso:order}
There is an order isomorphism 
\begin{equation} \label{Phi}
\Phi\colon (\varUpsilon([0,1],*), \prec)\to (\varUpsilon([0,1],\bullet), \prec)
\end{equation}
such that
\begin{equation} \label{Phi-def}
(a,b)\in\CI_i^*\iff\Phi(a,b)\in\CI_i^{\bullet}\quad (i\in\{P,\text{\L},M\}).
\end{equation}
\end{enumerate}
\end{thm}

\begin{proof}
\ref{t-norm-iso:cong}$\implies$\ref{t-norm-iso:order}: Let $\phi\colon([0,1],*)\to([0,1],\bullet)$
be an isomorphism of continuous t-norms. Define
\[\Phi\colon\varUpsilon([0,1],*)\to\varUpsilon([0,1],\bullet),\quad (a,b)\mapsto (\phi(a),\phi(b)).\]
Then $\Phi$ is clearly an order isomorphism because so is $\phi$. Moreover, \eqref{Phi-def} follows from the fact that for each $(a,b)\in\varUpsilon([0,1],*)$, the restriction $\phi|[a,b]\colon[a,b]\to[\phi(a),\phi(b)]$ is also an isomorphism of continuous t-norms.


\ref{t-norm-iso:order}$\implies$\ref{t-norm-iso:cong}: For each $(a,b)\in\varUpsilon([0,1],*)$, by \eqref{Phi-def} we find an isomorphism of continuous t-norms
\[\phi_{ab}\colon([a,b],*)\to([a_{\Phi},b_{\Phi}],\bullet),\]
where $(a_{\Phi},b_{\Phi}):=\Phi(a,b)$. Let $\psi_{ab}=\phi_{a,b}|(a,b)$, i.e., the restriction of $\phi_{ab}$ on $(a,b)$. Then
\[\psi=\bigcup\limits_{(a,b)\in\varUpsilon([0,1],*)}\psi_{ab}\colon\bigcup\varUpsilon([0,1],*)\to\bigcup\varUpsilon([0,1],\bullet)\]
is a well-defined and strictly increasing bijection, and it follows from Lemma \ref{t-norm-rep} that
\[\psi(x*y)=\psi(x)\bullet\psi(y)\]
for all $x,y\in\bigcup\varUpsilon([0,1],*)$. Since $\bigcup\varUpsilon([0,1],*)$ and $\bigcup\varUpsilon([0,1],\bullet)$ are both open dense subsets of $[0,1]$, the map
\[\phi\colon[0,1]\to[0,1],\quad\phi(x)=\sup\limits_{y\in[0,x]\cap(\bigcup\varUpsilon([0,1],*))}\psi(y)\]
is clearly an order isomorphism on $[0,1]$ that extends $\psi$, and it is also an isomorphism of monoids
\[\phi\colon([0,1],*,1)\to([0,1],\bullet,1)\]
by the continuity of $\phi$, $*$ and $\bullet$. Hence $([0,1],*)\cong_t([0,1],\bullet)$.
\end{proof}

\begin{exmp} \label{limit}
Continuous t-norms of the forms
\begin{equation} \label{1-1(n+1),1-1(n+2)}
\Big(\Big[1-\dfrac{1}{n+1}, 1-\dfrac{1}{n+2}\Big], *\Big)_{n\in \omega}\quad\text{and}\quad\Big(\Big[\dfrac{1}{n+2}, \dfrac{1}{n+1}\Big], \bullet\Big)_{n\in \omega}
\end{equation}
are not isomorphic, even if $\Big(\Big[1-\dfrac{1}{n+1}, 1-\dfrac{1}{n+2}\Big], *\Big)$ and $\Big(\Big[\dfrac{1}{n+2}, \dfrac{1}{n+1}\Big], \bullet\Big)$ are both isomorphic to $([0,1],*_P)$ (or, both isomorphic to $([0,1],*_{\text{\L}})$) for all $n\in\omega$. Indeed, we have
\[\CI_P^* =\Big\{\Big(1-\dfrac{1}{n+1}, 1-\dfrac{1}{n+2}\Big)\mathrel{\Big|} n\in\omega\Big\},\quad\CI_P^{\bullet}=\Big\{\Big(\dfrac{1}{n+2}, \dfrac{1}{n+1}\Big)\mathrel{\Big|} n\in\omega\Big\}\quad\text{and}\quad\CI_{\text{\L}}^*=\CI_M^*=\CI_{\text{\L}}^{\bullet}=\CI_M^{\bullet}=\varnothing.\] 
Thus, $\varUpsilon([0,1],*)=\CI_P^*$ has a $\prec$-minimum element, i.e., $\Big(0,\dfrac{1}{2}\Big)$, while $\varUpsilon([0,1],\bullet)=\CI_P^{\bullet}$ does not. Therefore, there is no order isomorphism from $\varUpsilon([0,1],*)$ to $\varUpsilon([0,1],\bullet)$; by Theorem \ref{t-norm-iso}, this means that the continuous t-norms \eqref{1-1(n+1),1-1(n+2)} cannot be isomorphic.
\end{exmp}

\begin{exmp} \label{Cantor}
Let $2=\{0,1\}$. Write $2^{<\omega}$ for the binary tree, i.e., the set of finite strings made up of $0$ and $1$. Recall that a family $J=\{J_u\}_{u\in 2^{<\omega}}$ of sets is called a \emph{Cantor System} (see, e.g., \cite[Definition 899]{Dasgupta2014})
if for each $u\in 2^{<\omega}$:
\begin{itemize}
\item $J_u=[a,b]$ for some $a<b$ in $\bbR$;
\item $J_{u^\smallfrown 0}$, $J_{u^\smallfrown 1} \subseteq J_u$, where $u^\smallfrown 0$ and $u^\smallfrown 1$ refer to the extensions of the finite string $u$ obtained by adding $0$ and $1$ at the end of $u$, respectively;
\item $J_{u^\smallfrown 0}\cap J_{u^\smallfrown 1}=\varnothing$;
\item $\lim\limits_{n\rightarrow\infty}|J_{b_0 b_1 \dots b_n}|=0$, where $b=b_0 b_1\dots b_n\dots\in 2^{\omega}$ is an arbitrary infinite binary sequence, and $|J_{b_0 b_1 \dots b_n}|$ refers to the length of the interval $J_{b_0 b_1 \dots b_n}$.
\end{itemize}
A set $A\subseteq\bbR$ is called a \emph{generalized Cantor set} if it is generated by some Cantor system $J=\{J_u\}_{u\in 2^{<\omega}}$; explicitly, 
\[x\in A\iff \exists b=b_0 b_1 \dots b_n\dots\in 2^{\omega}\colon x\in\bigcap\limits_{n\in\omega}J_{b_0 b_1 \dots b_n}.\] 
For example, the (classical) \emph{Cantor set} and the \emph{$\varepsilon$-Cantor set} (also \emph{Smith--Volterra--Cantor set}, see, e.g., \cite[Section 18]{Aliprantis1998} and \cite[Exercise 3.24]{Vallin2013}) are both generalized Cantor sets. 

For each generalized Cantor set $A$ generated by a Cantor system $J=\{J_u\}_{u\in 2^{<\omega}}$ with $J_{\varnothing}=[0,1]$, we denote by
\[\CK_A=\{(a_n,b_n)\}_{n\in\omega}\] 
the collection of the open intervals removed from $[0,1]$ during the construction of $A$. Note that the open intervals in $\CK_A$ are nonempty and pairwise disjoint. Define a continuous t-norm 
\begin{equation} \label{*A-def}
([0,1],*_A):= ([a_n,b_n],*_A)_{n\in\omega}
\end{equation}
with each $([a_n,b_n],*_A)$ being isomorphic to $([0,1], *_P)$ ($n\in \omega$). From the construction of \eqref{*A-def} it is easy to see that
\begin{equation} \label{*A-IP-IL-IM}
\CI_P^{*_A}=\CK_A,\quad\CI_{\text{\L}}^{*_A}=\varnothing\quad\text{and}\quad\bigcup\CI_M^{*_A}\subseteq A.
\end{equation}

If $A$ is the Cantor set or the $\varepsilon$-Cantor set, it is obvious that the associated Cantor system $J=\{J_u\}_{u\in 2^{<\omega}}$ has the following property:
\begin{enumerate}[label=(E)]
\item \label{E} For each $u\in 2^{<\omega}$, $J_u$ and $J_{u^\smallfrown 0}$ have the same left endpoint, while $J_u$ and $J_{u^\smallfrown 1}$ have the same right endpoint.
\end{enumerate}
Intuitively, the property \ref{E} means that during the construction of $A$, the open interval removed at each step is always located in the middle of the original closed interval (i.e., the left and right endpoints of the removed interval and the original interval are always different). Therefore, in this case, the order $\prec$ on $\CK_A$ (defined in the same way as in \eqref{CS-order}) has no endpoints; that is, $(\CK_A,\prec)$ has no minimum or maximum elements.

{\bf Fact 1.} Let $A$ and $B$ be generalized Cantor Sets generated by the Cantor systems $J= \{J_u\}_{u\in 2^{<\omega}}$ and $L=\{L_v\}_{v\in 2^{<\omega}}$, respectively, where $J_{\varnothing}=L_{\varnothing}=[0,1]$. If both $J$ and $L$ have the property \ref{E}, then the continuous t-norms $([0,1], *_A)$ and $([0,1], *_B)$ are isomorphic.

\emph{Proof of Fact 1.} Since $A$ and $B$ are generalized Cantor Sets and $J_{\varnothing}=L_{\varnothing}=[0,1]$, they are both nowhere dense (by \cite[Therorem 977]{Dasgupta2014}), which in conjunction with \eqref{*A-IP-IL-IM} leads to
\begin{align*}
\varUpsilon([0,1],*_A)&=\CI_P^{*_A}\cup\CI_{\text{\L}}^{*_A}\cup\CI_M^{*_A}= \CK_A \cup \varnothing \cup \varnothing = \CK_A,\\
\varUpsilon([0,1],*_B)&=\CI_P^{*_B}\cup\CI_{\text{\L}}^{*_B}\cup\CI_M^{*_B}= \CK_B \cup \varnothing \cup \varnothing = \CK_B.
\end{align*}
Since both $J$ and $L$ have the property \ref{E}, $(\varUpsilon([0,1],*_A),\prec)=(\CK_A,\prec)$ and $(\varUpsilon([0,1],*_B),\prec)=(\CK_B,\prec)$ are countable
dense orders without endpoints, which are necessarily isomorphic by \cite[Therorem 541]{Dasgupta2014}. Hence, Theorem \ref{t-norm-iso} guarantees that $([0,1], *_A)$ and $([0,1], *_B)$ are isomorphic. \qed

{\bf Fact 2.} Let $A$ and $B$ be generalized Cantor Sets generated by the Cantor systems $J= \{J_u\}_{u\in 2^{<\omega}}$ and $L=\{L_v\}_{v\in 2^{<\omega}}$, respectively, where $J_{\varnothing}=L_{\varnothing}=[0,1]$. If $J$ has the property \ref{E} but $L$ does not, 
then the continuous t-norms $([0,1], *_A)$ and $([0,1], *_B)$ are not isomorphic.

\emph{Proof of Fact 2.} Analogously to the above proof, we know that $(\varUpsilon([0,1],*_A),\prec)=(\CK_A,\prec)$ is a countable
dense order without endpoints, and $\varUpsilon([0,1],*_B)=\CK_B$.  By Theorem \ref{t-norm-iso}, it suffices to show that the orders $(\CK_A,\prec)$ and $(\CK_B,\prec)$ are not isomorphic when $L$ does not have the property \ref{E}. To this end, we show that $(\CK_B,\prec)$ either has an endpoint or is not dense.

Let $l(I)$ (resp. $r(I)$) denote the left (resp. right) endpoint of an interval $I$, respectively. Since $L$ does not have the property \ref{E}, there exists $v\in 2^{<\omega}$ such that $l(L_v)\neq l(L_{v^{\smallfrown} 0})$ or $r(L_v)\neq r(L_{v^{\smallfrown} 1})$. Suppose that $l(L_v)\neq l(L_{v^{\smallfrown} 0})$. Then there are two cases:
\begin{itemize}
\item $l(L_v)=0$. In this case, $(0, l(L_{v^{\smallfrown} 0})$ is the minimum element of $(\CK_B,\prec)$.
\item $l(L_{v^{\smallfrown} 0})>l(L_v)>0$. In this case, by the construction of $B$, there exists $(a_n, b_n)\in\CK_B$ such that $b_n=l(L_v)$. Thus, $(l(L_v),l(L_{v^{\smallfrown} 0}))$ is an immediate successor of $(a_n, b_n)$ in $(\CK_B,\prec)$, which means that $(\CK_B,\prec)$ is not dense.
\end{itemize}
If $r(L_v)\neq r(L_{v^{\smallfrown} 1})$, it can be proved similarly that $(\CK_B,\prec)$ either has a maximum element or is not dense. \qed
\end{exmp}

\section{Reducibility of equivalence relations} \label{Reducibility}

For the convenience of readers who may not be familiar with descriptive set theory, in this section we provide a brief review of basic facts about the reducibility of equivalence relations that will be used later.

A topological space is \emph{Polish} if it is separable and completely metrizable. The following properties of Polish spaces are well known (see, e.g., \cite[Proposition 3.3]{Kechris1995}): 

\begin{prop} \label{closed-product-Polish}
A closed subspace of a Polish space is Polish; a countable product of Polish spaces is Polish.
\end{prop}


Note that by \cite[Theorem 4.19]{Kechris1995}, the set $\CC$ of continuous functions from $[0,1]\times[0,1]$ to $[0,1]$ becomes a Polish space by taking the sup metric
\begin{equation} \label{sup-metric-def}
d_{\infty}(f_1,f_2)=\sup\{|f_1(x,y)-f_2(x,y)|\mid(x,y)\in[0,1]\times[0,1]\}.
\end{equation}
In what follows we always assume that $\CC$ is equipped with the topology induced by the metric $d_{\infty}$.

\begin{prop} \label{T-Polish}
The subset $\CT$ of $\CC$ consisting of continuous t-norms on $[0,1]$ is closed. Hence, as a subspace of $\CC$, $\CT$ is also a Polish space.
\end{prop}

\begin{proof}
It is straightforward to check that the limit of a sequence of continuous t-norms on $[0,1]$ in $\CC$ is also a continuous t-norm. Thus $\CT$ is a closed subspace of $\CC$, and hence it is Polish by Proposition \ref{closed-product-Polish}.
\end{proof}


Let $X$, $Y$ be Polish spaces, and let $E$, $F$ be equivalence relations on $X$, $Y$, respectively. We record some basic notions about Borel reducibility \cite{Kechris1995,Gao2008}:
\begin{enumerate}[label=(B\arabic*)]
\item \label{Borel:set} $A\subseteq X$ is \emph{Borel} if it belongs to the smallest $\sigma$-algebra on $X$ containing all open sets of $X$. 
\item \label{Borel:analytic} $A\subseteq X$ is \emph{analytic} (or $\mathbf{\Sigma}^1_1$) if there exists a Polish space $Z$ and a continuous function $f\colon Z\to X$ such that $f(Z)=A$. 
\item \label{Borel:complete-analytic} $A\subseteq X$ is \emph{complete analytic} if $A$ is analytic, and for each analytic subset $C$ of a zero-dimensional Polish space $Z$,  there exists a continuous function $f\colon Z\to X$ such that $C=f^{-1}(A)$.
\item \label{Borel:function} A function $f\colon X\to Y$ is \emph{Borel} if $f^{-1}(V)$ is Borel for any Borel (equivalently, open or closed) set $V$ of $Y$. In particular, if $Y$ has a countable subbasis $\{V_n\}_{n\in\omega}$, then $f$ is Borel provided that $f^{-1}(V_n)$ is Borel for all $n\in\omega$.
\item \label{Borel:reduction} A function $f\colon X\to Y$ is a \emph{reduction} from $(X,E)$ to $(Y,F)$ if
\[x_1Ex_2 \iff f(x_1)Ff(x_2)\]
for all $x_1$, $x_2\in X$.
\item \label{Borel:Borel-reducible} $(X,E)$ is  \emph{Borel reducible} to $(Y,F)$, denoted by
\[(X,E)\leq_B(Y,F),\]
if there exists a Borel function $f\colon X\to Y$ such that $f\colon(X,E)\to(Y,F)$ is a reduction.
\item \label{Borel:Borel-bireducible} $(X,E)$ and $(Y,F)$ are \emph{Borel bireducible} to each other, denoted by
\[(X,E)\sim_B(Y,F),\]
if $(X,E)\leq_B(Y,F)$ and $(Y,F)\leq_B(X,E)$.
\item \label{Borel:smooth} $(X,E)$ is \emph{smooth} (or \emph{concretely classifiable}) (see \cite[Definition 5.4.1]{Gao2008}), if
\[(X,E)\leq_B(2^\omega, \id(2^\omega)),\] 
where $2=\{0,1\}$ is equipped with the discrete topology, and $\id(2^\omega)$ is the identity relation on $2^\omega$; that is, if there exists a Borel function 
\begin{equation} \label{Borel:smooth:function}
f\colon X\to2^\omega\quad\text{with}\quad xEy\iff f(x)=f(y).
\end{equation}
\item \label{Borel:complete} $(X,E)$ is \emph{Borel complete} (cf. \cite[Definition 13.1.1]{Gao2008}), if it is Borel bireducible with the universal $S_{\infty}$-orbit equivalence relation.
\end{enumerate}

\begin{exmp} \label{Borel-complete-exmp}
The following equivalence relations are Borel complete:
\begin{itemize}
\item the isomorphism relation between countable graphs (see \cite[Theorem 13.1.2]{Gao2008});
\item the isomorphism relation between countable groups (see \cite[Theorem 13.4.1]{Gao2008});
\item the isomorphism relation between countable Boolean algebras \cite{Camerlo2001};
\item the homeomorphism relation between separable Boolean spaces (i.e., zero-dimensional compact metrizable spaces) \cite{Camerlo2001};
\item the isomorphism relation between commutative almost finite-dimensional $C^*$-algebra \cite{Camerlo2001};
\item the equivalence relation of isometry between Polish ultrametric spaces (see \cite[Theorem 4]{Gao2003});
\item the isomorphism relation between countable torsion-free abelian groups with domain $\omega$ \cite{Paolini2024}. 
\end{itemize}
\end{exmp}

The following lemma will be used later, which states that the property of being a complete analytic set is invariant under Borel bireducibility. The proof is elementary and is left to the reader.

\begin{lem} \label{Borel-bireducible-analytic}
Suppose that $(X,E)$ is Borel bireducible with $(Y,F)$, where $E$ and $F$ are equivalence relations on Polish spaces $X$ and $Y$, respectively. Then $E$ is a complete analytic subset of $X\times X$ (under the product topology) if, and only if, $F$ is a complete analytic subset of $Y\times Y$.
\end{lem}

Recall that a finite \emph{relational language} is a set $L=\{R_i\}_{i\in I}$, where $I$ is finite and each $R_i$ is assigned an element $n_i\in\omega$ $(i\in I)$, meaning that $R_i$ is an $n_i$-ary \emph{relation symbol}. Let $2=\{0,1\}$ be equipped with the discrete topology. We denote by 
\begin{equation} \label{ModL-product}
\Mod(L):=\prod\limits_{i\in I}2^{\omega^{n_i}}
\end{equation}
the set of \emph{$L$-structures} with the underlying set $\omega$, whose elements are of the form
\begin{equation} \label{RiS-def}
S=\{(R_i^S)_{i\in I}\mid\forall i\in I\colon R_i^S\subseteq\omega^{n_i}\}\in\prod\limits_{i\in I}2^{\omega^{n_i}}.
\end{equation}
As an immediate consequence of Lemma \ref{closed-product-Polish}, $\Mod(L)$ is a compact Polish space under the product topology. We say that $S_1,S_2\in\Mod(L)$ are \emph{isomorphic}, denoted by
\[S_1\cong_L S_2,\]
if there exists a bijection $\xi\colon\omega\to\omega$ such that
\begin{equation} \label{cong-ModL-def}
(k_1,\dots,k_{n_i})\in R^{S_1}_i\iff (\xi(k_1),\dots,\xi(k_{n_i}))\in R^{S_2}_i
\end{equation}
for all $i\in I$ and $(k_1,\dots,k_{n_i})\in\omega^{n_i}$. 

Let 
\begin{equation} \label{L1-def}
L_1=\{<, R_P, R_{\text{\L}}, R_M\},
\end{equation}
where $<$ is a binary relation symbol and $R_P, R_{\text{\L}}, R_M$ are unary relation symbols. By \eqref{ModL-product} and \eqref{RiS-def}, each $S\in\Mod(L_1)$ is a quadruple
\[(<^S,R_P^S,R_{\text{\L}}^S,R_M^S)\in 2^{\omega^2}\times 2^{\omega}\times 2^{\omega}\times 2^{\omega}.\]

Let $L_0=\{<\}$, where $<$ is a binary relation symbol; let
\begin{equation} \label{LO-def}
\LO:=\{<^S\mid\ <^S\ \text{is a (strict) linear order on}\ \omega \}\subseteq\Mod(L_0).
\end{equation}
It is easy to check that $\LO$ is a closed subspace of $\Mod(L_0)$, and thus a Polish space by Lemma \ref{closed-product-Polish}. For $<^{S_1},<^{S_2}\in\LO$, it is clear that $<^{S_1}\cong_{L_0}<^{S_2}$ if and only if $(\omega,<^{S_1})$ and $(\omega,<^{S_2})$ are order isomorphic. Considering the restriction of the equivalence relation $\cong_{L_0}$ on $\LO$, the following result is well known:

\begin{lem} \label{L1-leqB-LO} (See \cite[Theorem 3]{Friedman1989} and \cite[Theorem 4.10]{Hjorth2000})
$(\Mod(L_1),\cong_{L_1})\sim_B(\LO,\cong_{L_0})$.
\end{lem}

\section{The complexity of the isomorphism relation between continuous t-norms}  \label{Complexity}





Enumerate elements of $\bbQ\cap[0,1]$ in a non-repeating sequence $\{q_n\}_{n\in\omega}$. For each continuous t-norm $([0,1],*)$, define
\begin{equation} \label{f*-nxy*-def}
n^*_{xy}:=\begin{cases}
\min{\{n\in\omega\mid q_n \in(x, y)\}} & \text{if}\ (x,y)\in\CI_P^*\cup \CI_{\text{\L}}^*,\\
\min{\{n\in\omega\mid q_n \in[x, y]\}} & \text{if}\ (x,y)\in\CI_M^*,
\end{cases}
\end{equation}
which intuitively replaces each interval $(x,y)\in \varUpsilon([0,1], *)$ with a natural number $n^*_{xy}$. Let
\[S_*=(<^{S_*},R_P^{S_*},R_{\text{\L}}^{S_*},R_M^{S_*})\]
be the $L_1$-structure (see \eqref{L1-def}) given by
\begin{equation} \label{S*-def-<}
m<^{S_*}n \iff \exists (x,y),(x',y')\in\varUpsilon([0,1],*)\colon n^*_{xy} =m\with n^*_{x'y'} =n\with(x,y)\prec(x',y')
\end{equation}
and
\begin{equation} \label{S*-def-PLM}
R_i^{S_*}=\{n^*_{xy}\mid(x,y)\in\CI_i^*\}\quad (i\in\{P,\text{\L},M\}).
\end{equation}

\begin{lem} \label{RPS*-RLS*-RMS*}
Let $m,n\in\omega$.
\begin{enumerate}[label={\rm(\arabic*)}]
\item \label{RPS*-RLS*-RMS*:P}
    $n\in R_P^{S_*}$ if, and only if, 
    \begin{enumerate}[label={\rm(\alph*)}]
    \item \label{RPS*-RLS*-RMS*:P:qn} $q_n$ is non-idempotent,
    \item \label{RPS*-RLS*-RMS*:P:qnl} $\forall l\in\omega\setminus\{0\}$, $(q_n)_*^{(l)}:=\underbrace{q_n*q_n*\dots*q_n}_{l\ \text{times}}$ is non-idempotent, and
    \item \label{RPS*-RLS*-RMS*:P:qi} $\forall i<n$, $q_i * q_n = \min\{q_i, q_n\}$.
    \end{enumerate}
\item \label{RPS*-RLS*-RMS*:L}
    $n\in R_{\text{\L}}^{S_*}$ if, and only if,  
    \begin{enumerate}[label={\rm(\alph*)}]
    \item \label{RPS*-RLS*-RMS*:L:qn} $q_n$ is non-idempotent,
    \item \label{RPS*-RLS*-RMS*:L:qnl} $\exists l\in\omega\setminus\{0\}$, $(q_n)_*^{(l)}$ is idempotent, and
    \item \label{RPS*-RLS*-RMS*:L:qi} $\forall i<n$, $q_i * q_n = \min\{q_i, q_n\}$.
    \end{enumerate}
\item \label{RPS*-RLS*-RMS*:M}
    $n\in R_M^{S_*}$ if, and only if,  
    \begin{enumerate}[label={\rm(\alph*)}]
    \item \label{RPS*-RLS*-RMS*:M:qn} $q_n$ is idempotent,
    \item \label{RPS*-RLS*-RMS*:M:qk} $\exists l\in \omega \setminus\{n\}$, $\forall k\in \omega$, $q_k\in(q_l, q_n)\cup(q_n, q_l)\implies q_k$ is idempotent, and
    \item \label{RPS*-RLS*-RMS*:M:qj} $\forall i < n$, $\exists j\in \omega$, $q_j\in (q_i, q_n)\cup(q_n, q_i)$ and $q_j$ is non-idempotent.
    \end{enumerate}
\item \label{RPS*-RLS*-RMS*:<}
    $m<^{S_*}n$ if, and only if, $m,n\in R_P^{S_*}\cup R_{\text{\L}}^{S_*}\cup R_M^{S_*}$ and $q_m<q_n$.
\end{enumerate}
\end{lem}

\begin{proof}
\ref{RPS*-RLS*-RMS*:P} If $n\in R_P^{S_*}$, then $n=n_{xy}^*$ for some $(x,y)\in\CI_P^*$ (by \eqref{S*-def-PLM}); that is, $q_n\in(x,y)$ and $q_i\not\in (x,y)$ for all $i<n$ (by \eqref{f*-nxy*-def}). Hence, \ref{RPS*-RLS*-RMS*:P:qn} and \ref{RPS*-RLS*-RMS*:P:qi} follow from Lemma \ref{t-norm-rep}. Moreover, \ref{RPS*-RLS*-RMS*:P:qnl} holds because $([x,y],*)\cong_t([0,1],*_P)$ implies that every $q\in(x,y)$ is non-idempotent and non-nilpotent. 

Conversely, if $q_n$ is non-idempotent, then there exists $(x,y)\in\CI_{P}^*\cup\CI_{\text{\L}}^*$ such that $q_n\in (x,y)$ (by Lemma \ref{t-norm-rep}). Note that $q_n$ is a non-nilpotent element of $([x,y],*)$ by \ref{RPS*-RLS*-RMS*:P:qnl}; otherwise, if $(q_n)_*^{(l)}=x$ for some $l\in\omega\setminus\{0\}$, then $(q_n)_*^{(l)}$ is idempotent, which violates \ref{RPS*-RLS*-RMS*:P:qnl}. Thus it follows from Lemma \ref{*-iso-PL-condition:L} that $(x,y)\not\in\CI_{\text{\L}}^*$, and consequently $(x,y) \in \CI_{P}^*$. Finally, by \ref{RPS*-RLS*-RMS*:P:qi} and Lemma \ref{*-iso-PL-condition:min} we have $q_i\not\in (x,y)$ for all $i<n$. Hence $n=n_{xy}^*\in R_P^{S_*}$.

\ref{RPS*-RLS*-RMS*:L} The proof is mostly similar to \ref{RPS*-RLS*-RMS*:P}, where the only difference is \ref{RPS*-RLS*-RMS*:L:qnl}. For the necessity, \ref{RPS*-RLS*-RMS*:L:qnl} follows from the fact that $q_n$ is nilpotent (by Lemma \ref{*-iso-PL-condition:L}). Conversely, for the sufficiency, \ref{RPS*-RLS*-RMS*:L:qnl} implies that the interval $(x,y)$ in $\CI_{P}^*\cup\CI_{\text{\L}}^*$ containing $q_n$ cannot be in $\CI_{P}^*$.



\ref{RPS*-RLS*-RMS*:M} If $n\in R_M^{S_*}$, then $n=n_{xy}^*$ for some $(x,y)\in\CI_M^*$ (by \eqref{S*-def-PLM}); that is, $q_n\in[x,y]$ and $q_i\not\in[x,y]$ for all $i<n$ (by \eqref{f*-nxy*-def}). Thus $q_n$ must be idempotent; because by the continuity of $*$, the endpoints of intervals in $\CI_M^*$ are necessarily idempotent. 

For \ref{RPS*-RLS*-RMS*:M:qj}, let $i<n$. Since $q_i\not\in[x,y]$, there exist non-idempotent elements in $(q_i, q_n)\cup(q_n, q_i)$ (by the definition of $\CI_M^*$). Thus, by the density of $\{q_j\}_{j\in\omega}$ in $[0,1]$ we find $j\in\omega$ such that $q_j\in (q_i, q_n)\cup(q_n, q_i)$ and $q_j$ is non-idempotent. 

For \ref{RPS*-RLS*-RMS*:M:qk}, since every $q\in[x,y]$ is idempotent, by the density of $\{q_j\}_{j\in\omega}$ in $[0,1]$ we may choose $l\in \omega \setminus\{n\}$ such that $q_l\in[x,q_n)\cup(q_n,y]$. Then $q_k$ is idempotent whenever $q_k\in(q_l, q_n)\cup(q_n, q_l)$. 

Conversely, from \ref{RPS*-RLS*-RMS*:M:qn} and \ref{RPS*-RLS*-RMS*:M:qk} we see that every $q\in\bbQ\cap((q_l, q_n)\cup(q_n, q_l))$ is idempotent, and thus, by Lemma \ref{t-norm-rep},  $(q_l, q_n)\cup(q_n, q_l)$ is a nonempty interval of idempotent elements. Therefore, there exists $(x,y)\in\CI_{M}^*$ such that $[q_l,q_n]\cup[q_n,q_l]\subseteq[x,y]$. Finally, \ref{RPS*-RLS*-RMS*:M:qj} guarantees $n_{xy}^*=n$, because there exists a non-idempotent element $q_j\in(q_i,q_n)\cup(q_n,q_i)$ for all $i<n$, which means that $q_i\not\in[x,y]$. Thus $n\in R_M^{S_*}$.

\ref{RPS*-RLS*-RMS*:<} If $m<^{S_*}n$, by \eqref{S*-def-<} we find $(x,y),(x',y')\in\varUpsilon([0,1],*)$ such that $n^*_{xy}=m$, $n^*_{x'y'} =n$ and $(x,y)\prec(x',y')$. Thus $m,n\in R_P^{S_*}\cup R_{\text{\L}}^{S_*}\cup R_M^{S_*}$ and $q_m\in(x,y)$, $q_n\in(x',y')$ (by \eqref{f*-nxy*-def} and \eqref{S*-def-PLM}), and consequently $q_m<q_n$ (by \eqref{CS-order}).

Conversely, if $m,n\in R_P^{S_*}\cup R_{\text{\L}}^{S_*}\cup R_M^{S_*}$ and $q_m<q_n$, then there exist $(x,y),(x',y')\in\varUpsilon([0,1],*)$ such that $n^*_{xy}=m$, $n^*_{x'y'} =n$ and $(x,y)\prec(x',y')$ (by \eqref{f*-nxy*-def} and \eqref{S*-def-PLM}), which means that $m<^{S_*}n$ (by \eqref{S*-def-<}).
%
%
%
\end{proof}

\begin{prop} \label{leqB-ModL1}
$(\CT,\cong_t)\leq_B(\Mod(L_1),\cong_{L_1})$.
\end{prop}

\begin{proof}
We show that
\[\Theta\colon(\CT,\cong_t)\to(\Mod(L_1),\cong_{L_1}),\quad\Theta([0,1], *)=S_*\]
is a Borel reduction. 

{\bf Step 1.} $\Theta$ is a reduction. Given continuous t-norms $([0,1],*)$ and $([0,1],\bullet)$, we verify that
\[([0,1],*)\cong_t([0,1],\bullet)\iff S_*\cong_{L_1}S_{\bullet}.\]

On one hand, if $([0,1],*)\cong_t([0,1],\bullet)$, then by Theorem \ref{t-norm-iso}, there is an order isomorphism \eqref{Phi} satisfying \eqref{Phi-def}. Write $(x_{\Phi},y_{\Phi}):=\Phi(x,y)\in\varUpsilon([0,1],\bullet)$ for each $(x,y)\in\varUpsilon([0,1],*)$. Note that both
\begin{equation} \label{omega-minus-f*}
\omega\setminus\{n_{xy}^*\mid(x,y)\in\varUpsilon([0,1],*)\}\quad\text{and}\quad\omega\setminus \{n_{x'y'}^{\bullet}\mid(x',y')\in\varUpsilon([0,1],\bullet)\}=\omega\setminus \{n^{\bullet}_{x_{\Phi}y_{\Phi}}\mid(x,y)\in\varUpsilon([0,1],*)\}
\end{equation}
are countably infinite sets, because each interval $(x,y)\in\varUpsilon([0,1],*)$ and the corresponding $(x_{\Phi},y_{\Phi})\in\varUpsilon([0,1],\bullet)$ have infinitely many elements of the sequence $\{q_n\}_{n\in\omega}$ (recall that $\{q_n\}_{n\in\omega}=\bbQ\cap[0,1]$), but in \eqref{omega-minus-f*} only one element is removed from each $(x,y)$ and $(x_{\Phi},y_{\Phi})$ (i.e., $n^*_{xy}$ and $n^{\bullet}_{x_{\Phi}y_{\Phi}}$), respectively. Thus, there exists a bijection 
\[\eta\colon\omega\setminus\{n_{xy}^*\mid(x,y)\in\varUpsilon([0,1],*)\}\to\omega\setminus \{n^{\bullet}_{x_{\Phi}y_{\Phi}}\mid(x,y)\in\varUpsilon([0,1],*)\}.\]
Define
\[\xi\colon\omega\to\omega,\quad\xi(n)=\begin{cases}
n^{\bullet}_{x_{\Phi}y_{\Phi}} & \text{if}\ n=n^*_{xy}\ \text{for some}\ (x,y)\in\varUpsilon([0,1],*),\\
\eta(n) & \text{else}.
\end{cases}\]
Then $\xi$ is a bijection, and it follows immediately from the definition of $S_*$ and $S_{\bullet}$ (see \eqref{S*-def-<} and \eqref{S*-def-PLM}) that $\xi$ satisfies \eqref{cong-ModL-def}. Hence $S_*\cong_{L_1}S_{\bullet}$.

On the other hand, if $S_*\cong_{L_1}S_{\bullet}$, then there exists a bijection $\xi\colon\omega\to\omega$ satisfying \eqref{cong-ModL-def}; that is,
\begin{equation} \label{RiS*-iff-RiSb}
n\in R_i^{S_*}\iff\xi(n)\in R_i^{S_{\bullet}}\quad (i\in\{P,\text{\L},M\}).
\end{equation}
Thus, for each $i\in\{P,\text{\L},M\}$ and $(x,y)\in\CI_i^*$, we have $n^*_{xy}\in R_i^{S_*}$, and consequently $\xi(n^*_{xy})\in R_i^{S_{\bullet}}$, which means that there exists $(x_{\Phi},y_{\Phi})\in\CI_i^{\bullet}$ such that $\xi(n^*_{xy})=n^{\bullet}_{x_{\Phi}y_{\Phi}}$ (by \eqref{S*-def-PLM}). Define
\[\Phi\colon\varUpsilon([0,1],*)\to\varUpsilon([0,1],\bullet),\quad \Phi(x,y)=(x_{\Phi},y_{\Phi}).\]
Then $\Phi$ is an order isomorphism by \eqref{S*-def-<}, and satisfies \eqref{Phi-def} by \eqref{S*-def-PLM} and \eqref{RiS*-iff-RiSb}. Hence $([0,1],*)\cong_t([0,1],\bullet)$ by Theorem \ref{t-norm-iso}.

{\bf Step 2.} $\Theta\colon\CT\to\Mod(L_1)$ is a Borel function. Fix a bijection
\[\omega\to\omega\times\omega,\quad n\mapsto(\iota(n),\kappa(n)).\]
Note that
\[\CD=\{U_{i,n}^k\mid i\in\{<,P,\text{\L},M\},\ n\in\omega,\ k\in\{0,1\}\}\]
is a countable subbasis of $\Mod(L_1)=2^{\omega^2}\times 2^{\omega}\times 2^{\omega}\times 2^{\omega}$, where 
\begin{itemize}
\item $U_{<,n}^0=\{(<^S,R_P^S,R_{\text{\L}}^S,R_M^S)\in\Mod(L_1)\mid\iota(n)\not<^S\kappa(n)\}$,\quad $U_{<,n}^1=\{(<^S,R_P^S,R_{\text{\L}}^S,R_M^S)\in\Mod(L_1)\mid\iota(n)<^S\kappa(n)\}$,
\item $U_{P,n}^0=\{(<^S,R_P^S,R_{\text{\L}}^S,R_M^S)\in\Mod(L_1)\mid n\not\in R_P^S\}$,\quad $U_{P,n}^1=\{(<^S,R_P^S,R_{\text{\L}}^S,R_M^S)\in\Mod(L_1)\mid n\in R_P^S\}$,
\item $U_{\text{\L},n}^0=\{(<^S,R_P^S,R_{\text{\L}}^S,R_M^S)\in\Mod(L_1)\mid n\not\in R_{\text{\L}}^S\}$,\quad $U_{\text{\L},n}^1=\{(<^S,R_P^S,R_{\text{\L}}^S,R_M^S)\in\Mod(L_1)\mid n\in R_{\text{\L}}^S\}$,
\item $U_{M,n}^0=\{(<^S,R_P^S,R_{\text{\L}}^S,R_M^S)\in\Mod(L_1)\mid n\not\in R_M^S\}$,\quad $U_{M,n}^1=\{(<^S,R_P^S,R_{\text{\L}}^S,R_M^S)\in\Mod(L_1)\mid n\in R_M^S\}$.
\end{itemize}
It suffices to check that $\Theta^{-1}(U_{i,n}^k)$ is Borel for all $U_{i,n}^k\in\CD$. To this end, we consider the following sets for $q\in [0,1]$ and $m,n,l\in \omega$:
\begin{itemize}
\item $\CV(q) = \{ ([0,1], *) \in \CT\mid q \text{ is an idempotent element of } *\}$;
\item $\CU(m,n) = \{ ([0,1], *) \in \CT\mid q_m*q_n = \min\{q_m, q_n\}\}$;
\item $\CV(l,n) = \CV((q_n)_*^{(l)})$;
\item $\CW(m,n) = \{([0,1], *) \in \CT \mid q_m < q_n\}$.
\end{itemize}
It is easy to see that $\CV(q)$, $\CU(m,n)$ and $\CV(l,n)$ are all closed subsets of the metric space $(\CT,d_{\infty})$ (cf. \eqref{sup-metric-def}), and $\CW(m,n)$ is either $\CT$ or $\varnothing$. Therefore, the conclusion is an immediate consequence of the following expressions obtained by Lemma \ref{RPS*-RLS*-RMS*}:
\begin{align*}
\Theta^{-1}(U^1_{P,n})={}&(\CT\setminus\CV(q_n))\cap\Big(\bigcap_{l\in\omega\setminus\{0\}}(\CT\setminus\CV(l, n)) \Big) \cap \Big( \bigcap_{i<n} \CU(i, n) \Big),\\
\Theta^{-1}(U^1_{\text{\L},n})={}&(\CT\setminus\CV(q_n)) \cap \Big( \bigcup_{l\in\omega\setminus\{0\}} \CV(l, n) \Big) \cap \Big( \bigcap_{i<n} \CU(i, n) \Big),\\
\Theta^{-1}(U^1_{M,n})={}&\CV(q_n)\cap 
\Big( \bigcup_{l\in \omega\setminus\{n\}} \bigcap\{ \CV(q)\mid q \in \bbQ\cap((q_l, q_n)\cup(q_n, q_l))\} \Big)\\
&\cap \Big( \bigcap_{i<n}\bigcup \{\CT\setminus\CV(q)\mid q\in \bbQ\cap((q_{i}, q_n)\cup(q_n, q_{i}))\}\Big),\\
\Theta^{-1}(U^1_{<,n})={}&\Big( \bigcup_{i\in\{P,\text{\L},M\}} \Theta^{-1}(U^1_{i,\iota(n)}) \Big) \cap \Big( \bigcup_{i\in\{P,\text{\L},M\}} \Theta^{-1}(U^1_{i,\kappa(n)}) \Big) \cap \CW(\iota(n), \kappa(n)),\\
\Theta^{-1}(U^0_{i,n})={}&\Theta^{-1}(\Mod(L_1)\setminus U^1_{i,n}) =\CT\setminus\Theta^{-1}(U^1_{i,n})\quad (i\in\{<, P,\text{\L},M\}). 
\qedhere
\end{align*}
\end{proof}

Since $\cong_{L_1}$ is analytic, the following corollary is a direct consequence of Proposition \ref{leqB-ModL1} and \cite[Proposition 14.4]{Kechris1995}:

\begin{cor} \label{cong-t-analytic}
$\cong_t$ is an analytic equivalence relation on $\CT$.
\end{cor}

Since the Borel reducibility $\leq_B$ is transitive, from Lemma \ref{L1-leqB-LO} and Proposition \ref{leqB-ModL1} we immediately have:

\begin{prop} \label{leqB-congLO}
$(\CT,\cong_t)\leq_B(\LO,\cong_{L_0})$.
\end{prop}

Conversely, we have the following proposition:

\begin{prop} \label{L-leqC}
$(\LO,\cong_{L_0})\leq_B(\CT,\cong_t)$.
\end{prop}

\begin{proof}
{\bf Step 1.} For each (strict) linear order $\lessdot$ on $\omega$, we define a sequence of disjoint open intervals $\{I_n^{\lessdot}=(a_n^{\lessdot},b_n^{\lessdot})\}_{n\in\omega}$ such that the following properties hold for all $m,n\in\omega$:
\begin{enumerate}[label={\rm(\alph*)}]
\item \label{I-condition:01}
$0<a_n^{\lessdot}<b_n^{\lessdot}<1$;
\item \label{I-condition:endpoint}
$b_m^{\lessdot}<a_n^{\lessdot}$ whenever $m\lessdot n$;
\item \label{I-condition:length}
$|I_n|=\dfrac{1}{3^{n+1}}$, where $|I_n|$ refers to the length of $I_n$;
\item \label{I-condition:distance}
$\min\{a_n^{\lessdot},\ 1-b_n^{\lessdot},\ a_n^{\lessdot}-b_m^{\lessdot}\}\geq\dfrac{1}{3^{n+1}}$ whenever $m\lessdot n$.
\end{enumerate}
We proceed by induction. First, let 
\[I_0=(a_0^{\lessdot},b_0^{\lessdot})=\Big(\dfrac{1}{3},\dfrac{2}{3}\Big).\] 
Second, we consider
\[x_1=\begin{cases}
0 & \text{if}\ 1\lessdot 0\\
b_0^{\lessdot} & \text{if}\ 0\lessdot 1
\end{cases}
\quad\text{and}\quad
y_1=\begin{cases}
a_0^{\lessdot} & \text{if}\ 1\lessdot 0\\
1 & \text{if}\ 0\lessdot 1,
\end{cases}\]
and define $I_1^{\lessdot}=(a_1^{\lessdot}, b_1^{\lessdot})$ with 
\[a_1^{\lessdot}=\dfrac{1}{2}\Big(x_1+y_1-\dfrac{1}{9}\Big)\quad\text{and}\quad b_1^{\lessdot}=\dfrac{1}{2}\Big(x_1+y_1+\dfrac{1}{9}\Big),\] 
which clearly satisfies \ref{I-condition:01}--\ref{I-condition:distance}. In general, suppose that we have defined $I_k^{\lessdot}$ $(0\leq k\leq n-1)$ ($n\in\omega\setminus\{0\}$) satisfying \ref{I-condition:01}--\ref{I-condition:distance}. Then we consider
\[x_n=\max(\{0\}\cup\{b_k^{\lessdot}\mid k\lessdot n \with k< n\})\quad\text{and}\quad y_n=\min(\{a_k^{\lessdot}\mid n\lessdot k \with k< n\}\cup\{1\}).\] 
By \ref{I-condition:distance}, we can see that 
\[y_n-x_n\geq \dfrac{1}{3^n}.\]
Then we define $I_n^{\lessdot}=(a_n^{\lessdot}, b_n^{\lessdot})$ with 
\[a_n^{\lessdot}=\dfrac{1}{2}\Big(x_n+y_n-\dfrac{1}{3^{n+1}}\Big)\quad\text{and}\quad b_n^{\lessdot}=\dfrac{1}{2}\Big(x_n+y_n+\dfrac{1}{3^{n+1}}\Big).\] 
It is straightforward to check that $I_n$ also satisfies \ref{I-condition:01}--\ref{I-condition:distance}. Thus, we have $I_n^{\lessdot}=(a_n^{\lessdot},b_n^{\lessdot})$ defined for all $n\in\omega$.

{\bf Step 2.} For each $\lessdot\in\LO$, by \eqref{t-norm-def-by-rep} we may define a continuous t-norm 
\begin{equation} \label{[0,1]*<-def}
([0,1],*^{\lessdot})=([a_n^{\lessdot},b_n^{\lessdot}],*^{\lessdot})_{n\in\omega}
\end{equation}
where
\[x*^{\lessdot}y=a_n^{\lessdot}+\dfrac{(x-a_n^{\lessdot})(y-a_n^{\lessdot})}{b_n^{\lessdot}-a_n^{\lessdot}}\]
for all $x,y\in[a_n^{\lessdot},b_n^{\lessdot}]$. For any $n\in\omega$, it is easy to see that $([a_n^{\lessdot},b_n^{\lessdot}],*^{\lessdot})$ is isomorphic to the product t-norm $([0,1],*_P)$, and it is completely determined by the interval $I_n^{\lessdot}=(a_n^{\lessdot},b_n^{\lessdot})$; that is, for any $*^{\lessdot_1},*^{\lessdot_2}\in\LO$, 
\begin{equation} \label{an1-an2-bn1-bn2}
I_n^{\lessdot_1}=I_n^{\lessdot_2}\implies([a_n^{\lessdot_1},b_n^{\lessdot_1}],*^{\lessdot_1})=([a_n^{\lessdot_2},b_n^{\lessdot_2}],*^{\lessdot_2}).
\end{equation}
We show that the map
\[\theta\colon\LO\to\CT,\quad \lessdot\mapsto([0,1],*^{\lessdot})\]
is continuous, and thus Borel. To this end, for each $\lessdot_0\in\LO$ and $\varepsilon>0$, we have to find an open neighborhood $V$ of $\lessdot_0$ in $\LO$ such that 
\[\theta(V)\subseteq B(([0,1],*^{\lessdot_0}),\varepsilon),\] 
where the latter is an open ball in the metric space $(\CT,d_{\infty})$ (cf. \eqref{sup-metric-def}). Indeed, it follows from \ref{I-condition:length} that there exists $N\in\omega\setminus\{0\}$ such that 
\begin{equation} \label{sum-In-ep}
\sum\limits_{n\geq N}|I_n^{\lessdot}|<\dfrac{\varepsilon}{2}
\end{equation}
for all $\lessdot \in \LO$. By \eqref{ModL-product} and \eqref{LO-def},
\[V=\{\lessdot\in\LO\mid\lessdot |_{N\times N} = \lessdot_0 |_{N \times N}\}=\{\lessdot\in\Mod(L_0)=2^{\omega\times\omega}\mid\lessdot |_{N\times N} = \lessdot_0 |_{N \times N}\}\cap\LO\]
is open in $\LO$, where $N=\{n\in\omega\mid 0\leq n\leq N-1\}$ and $\lessdot |_{N\times N}=\lessdot\cap(N\times N)$. Since it is clear that $\lessdot_0\in V$, it remains to show that 
\begin{equation} \label{dinfty-*<-*<0-ep}
d_{\infty}(*^{\lessdot},*^{\lessdot_0})=\sup\{|x*^{\lessdot}y-x*^{\lessdot_0}y|\mid x,y\in[0,1]\}<\varepsilon
\end{equation}
for all $\lessdot\in V$. Indeed, by the constructions of $V$, $\{I_n^{\lessdot}\}_{n\in\omega}$ and $\{I_n^{\lessdot_0}\}_{n\in\omega}$ we have
\begin{equation} \label{In-n<N-equal}
\bigcup\limits_{n<N}I_n^{\lessdot}=\bigcup\limits_{n<N}I_n^{\lessdot_0}.
\end{equation}
Thus, if $x*^{\lessdot}y\neq x*^{\lessdot_0}y$, then either $x,y\in I_m^{\lessdot}$ or $x,y\in I_m^{\lessdot_0}$ for some $m\geq N$. Otherwise,
\begin{itemize}
\item if there exists no $n\in\omega$ such that $x,y\in I_n^{\lessdot}$ or $x,y\in I_n^{\lessdot_0}$, then $x*^{\lessdot}y=\min\{x,y\}=x*^{\lessdot_0}y$ (by Lemma \ref{t-norm-rep});
\item if $x,y\in I_n^{\lessdot}=I_n^{\lessdot_0}$ for some $n<N$ (by \eqref{In-n<N-equal}), then $x*^{\lessdot}y=x*^{\lessdot_0}y$ (by \eqref{an1-an2-bn1-bn2}).
\end{itemize}
So, there are three cases:
\begin{itemize}
\item $x,y\in I_m^{\lessdot}$ for some $m\geq N$, but $(x,y)\not\in I_n^{\lessdot_0}$ for all $n\geq N$. In this case, we have
\[x*^{\lessdot}y\in(a_m,b_m)\quad\text{and}\quad x*^{\lessdot_0}y=\min\{x,y\}\in(a_m,b_m),\]
and it follows from \eqref{sum-In-ep} that
\begin{equation} \label{x*y-x*0y<ep2}
|x*^{\lessdot}y- x*^{\lessdot_0}y|\leq|I_m^{\lessdot}|<\dfrac{\varepsilon}{2}.
\end{equation}

\item $x,y\in I_m^{\lessdot_0}$ for some $m\geq N$, but $(x,y)\not\in I_n^{\lessdot}$ for all $n\geq N$. In this case, we may prove \eqref{x*y-x*0y<ep2} analogously to the first case.
\item $x,y\in I_m^{\lessdot}\cap I_n^{\lessdot_0}$ for some $m,n\geq N$. In this case, by \eqref{sum-In-ep} we also have
\[|x*^{\lessdot}y- x*^{\lessdot_0}y|\leq|I_m^{\lessdot_0}|+|I_n^{\lessdot}| <\dfrac{\varepsilon}{2}.\]
\end{itemize}
Therefore, the desired inequality \eqref{dinfty-*<-*<0-ep} is obtained.

{\bf Step 3.} $\theta\colon(\LO,\cong_{L_0})\to(\CT,\cong_t)$ is a reduction. Indeed, for any $\lessdot_1,\lessdot_2\in\LO$, by the construction \eqref{[0,1]*<-def} we immediately see that there exists an order isomorphism 
\[f\colon(\omega,\lessdot_1)\to(\omega,\lessdot_2)\]
if, and only if, there exists an order isomorphism 
\[\Phi\colon\varUpsilon([0,1],*^{\lessdot_1})\to\varUpsilon([0,1],*^{\lessdot_2})\]
satisfying \eqref{Phi-def}. Therefore, it follows from Theorem \ref{t-norm-iso} that 
\[\lessdot_1\cong_{L_0}\lessdot_2\iff([0,1],*^{\lessdot_1})\cong_t([0,1],*^{\lessdot_2}),\] which completes the proof.
\end{proof}

Our main result now arises from the combinations of Propositions \ref{leqB-congLO} and \ref{L-leqC}:

\begin{thm} \label{main}
$(\CT,\cong_t)\sim_B(\LO,\cong_{L_0})$.
\end{thm}

It is well known that the restriction of the isomorphism relation $\cong_{L_0}$ on $\LO$, i.e., 
$\cong_{L_0}\cap (\LO\times\LO)$, is a complete analytic subset of $\LO\times\LO$; see, e.g., \cite[Theorems 3 and 4]{Friedman1989}. Therefore, the following corollary is an immediate consequence of Lemma \ref{Borel-bireducible-analytic} and Theorem \ref{main}:

\begin{cor} \label{cong-t-complete-analytic}
The isomorphism relation $\cong_t$ is a complete analytic subset of $\CT\times\CT$.
\end{cor}

Note that for a smooth equivalence relation $E$ on a Polish space $X$ (see \ref{Borel:smooth}) equipped with a Borel function \eqref{Borel:smooth:function}, 
\[E=(f\times f)^{-1}(\id(2^\omega))\]
is necessarily a Borel subset of $X\times X$ under the product topology, because $\id(2^\omega)$ is closed in $2^\omega\times 2^\omega$ and it is easy to see that $f\times f$ is Borel. However, complete analytic subsets cannot be Borel (as a direct consequence of \cite[Corollary 26.2]{Kechris1995}), and thus Corollary \ref{cong-t-complete-analytic} implies the following:


\begin{cor} \label{congt-not-smooth}
The equivalence relation $\cong_t$ is not smooth.
\end{cor}

Since it is well known that $(\LO,\cong_{L_0})$ is Borel complete (see \ref{Borel:complete}, \cite[Theorem 3]{Friedman1989} and \cite[Theorem 13.3.2]{Gao2008}), we conclude from Theorem \ref{main} that:

\begin{cor} \label{congt-Borel-complete}
$(\CT,\cong_t)$ is Borel complete. Therefore, $(\CT,\cong_t)$ is Borel bireducible with each equivalence relation listed in Example \ref{Borel-complete-exmp}.
\end{cor}

\section*{Acknowledgment}

This work is supported by the National Natural Science Foundation of China (No. 12071319), the Science and Technology Department of Sichuan Province (No. 2023NSFSC1285) and the Fundamental Research Funds for the Central Universities (No. 2021SCUNL202). We thank the anonymous referee for several helpful remarks.


\end{document}